\documentclass[12pt]{article}
\usepackage{amsmath}
\usepackage{amsthm,amssymb}
\usepackage{graphicx}
\usepackage{enumerate}
\usepackage{apacite}
\usepackage{natbib}
\usepackage{url} 
\newtheorem{thrm}{Theorem}[section]
\usepackage{color}
\usepackage{algcompatible}
\usepackage{algorithm}
\usepackage{algpseudocode}
\usepackage{hyperref}
\newcommand{\blind}{0}

\newcommand{\ee}{{\mathbf{e}}}

\newcommand{\xx}{{\mathbf{x}}}
\newcommand{\yy}{{\mathbf{y}}}
\newcommand{\zz}{{\mathbf{z}}}

\newcommand{\bA}{{\mathbf{A}}}

\newcommand{\DD}{{\mathbf{D}}}

\newcommand{\II}{{\mathbf{I}}}

\newcommand{\LL}{{\mathbf{L}}}
\newcommand{\MM}{{\mathbf{M}}}

\newcommand{\WW}{{\mathbf{W}}}
\newcommand{\XX}{{\mathbf{X}}}

\newcommand{\ZZ}{{\mathbf{Z}}}

\newcommand{\be}{{\mbox{\boldmath $\beta$}}}

\newcommand{\zer}{{\mbox{\boldmath $0$}}}

\newcommand{\E}{\operatorname{E}}

\renewcommand{\mathbf}{\boldsymbol}

\newcommand{\trace}{\operatorname{trace}}

\begin{document}

\def\spacingset#1{\renewcommand{\baselinestretch}%
{#1}\small\normalsize} \spacingset{1}


\if0\blind
{
  \title{\bf On the effect of noise on fitting linear regression models}
  \author{Insha Ullah and A.H. Welsh\thanks{
    The authors gratefully acknowledge funding support from Data61, CSIRO. }\hspace{.2cm}\\
    Research School of Finance, Actuarial Studies and Statistics\\ Australian National University}
  \maketitle
} \fi


\if1\blind
{
  \bigskip
  \bigskip
  \bigskip
  \begin{center}
    {\LARGE\bf 
    On the impact of noise on fitting linear regression models
    }
\end{center}
  \medskip
} \fi

\bigskip
\begin{abstract}
In this study, we explore the effects of including noise predictors and noise observations when fitting linear regression models. We present empirical and theoretical results that show that double descent occurs in both cases, albeit with contradictory implications: the implication for noise predictors is that complex models are often better than simple ones, while the implication for noise observations is that simple models are often better than complex ones.  We resolve this contradiction by showing that it is not the model complexity but rather the implicit shrinkage by the inclusion of noise in the model that drives the double descent.  Specifically, we show how noise predictors or observations shrink the estimators of the regression coefficients and make the test error asymptote, and then how the asymptotes of the test error and the ``condition number anomaly'' ensure that double descent occurs.   
We also show that including noise observations in the model makes the (usually unbiased) ordinary least squares estimator biased and indicates that the ridge regression estimator may need a negative ridge parameter to avoid over-shrinkage.

\end{abstract}

\vspace{3 cm}
\noindent%
{\it Keywords:}  Bias-variance tradeoff, Overparameterization, Multiple-descent, Negative ridge parameter; Regularization; Shrinkage; Sparsity.
\vfill

\newpage
\spacingset{1.9} 
\section{Introduction}
\label{sec:intro}

Statistical learning, particularly with large, high-dimensional datasets, currently rests on two opposing ideas, namely sparsity and overfitting.  Sparsity is the idea that a small (unknown) number of predictors are related to the response, and the remaining predictors are irrelevant noise variables that should be identified and discarded. Overfitting is the idea that very complex models (i.e., models with a very large number of predictors and/or parameters) can perform better than simpler models in prediction and classification. In this paper, we explore the effect of overfitting linear regression models when sparsity holds.

We consider the prediction task in which we use predictors $\mathbf{x'}$ to predict the response or label $y'$ in an observation $(y',\mathbf{x'})$ drawn from an unknown distribution $\mathcal{D}$.  Suppose we have an independent training dataset consisting of $n$ observations $(y_1,\mathbf{x}_1)\ldots, (y_n, \mathbf{x}_n) \in \mathbb{R} \times\mathbb{R}^{d}$ drawn independently from $\mathcal{D}$.  We use the training data to fit a working linear regression model to the test data by calculating $\hat{\boldsymbol{\beta}}\in \mathbb{R}^{ d}$, and then predict $y'$ by $\mathbf{x'}^T\hat{\boldsymbol{\beta}}$.    We take $\hat{\boldsymbol{\beta}}$ to be the least squares estimator when $d < n$ (the underparameterized regime) and the minimum norm least squares estimator when $d \ge n$ (the overparameterized regime).  Suppose we also have a test data set $(y'_1,\mathbf{x}'_1)\ldots, (y'_{n'},\mathbf{x}'_{n'})$ of size $n'$ drawn independently from $\mathcal{D}$ that is independent of the training data. Then we compare different models (e.g. using different elements of $\mathbf{x'}$) through their test error  $n'^{-1}\sum_{i=1}^{n'}(y'_i - \mathbf{x}'^{T}_i\hat{\boldsymbol{\beta}})^2$ and base the prediction on the model with the minimum test error.  

Under sparsity, only $d_0$ of the available predictors are related to the response, and the remaining variables are irrelevant noise variables. When $d_0 <<n$, the test error over the underparameterized regime ($d<n$) is typically roughly U-shaped; underfitting ($d<d_0$) by excluding variables related to the response increases the bias in the test error while overfitting by including irrelevant noise variables in the model ($d>d_0$) increases the variance in the test error. Selecting the model with minimum test error is intended to produce a model of dimension $d=d_0$ that includes only the relevant predictors by balancing the effects of under- and over-fitting, referred to as the bias-variance trade-off. There are many different methods for model selection and a vast literature on these methods.

Recent work exploring increasing $d$ past the interpolation point ($d=n$) into the overparameterized regime ($d>n$) has found that the test error decreases again, exhibiting a second descent with increasing model complexity \citep{belkin2019reconciling, geiger2019jamming, belkin2020two, bartlett2020benign, kobak2020optimal, liang2020just, liang2020multiple, muthukumar2020harmless, nakkiran2021deep, zhang2021understanding, hastie2022surprises, deng2022model, rocks2022memorizing}. 
 The second minimum (in the overparameterized regime) can be smaller than the first minimum (in the underparameterized regime), implying that large models with $d > n$ may give better predictions. This is called ``double descent'' by \cite{belkin2019reconciling}; ``triple'' \citep{nakkiran2020optimal} and even ``multiple descent'' \citep{liang2020multiple} have been demonstrated.  
 Double descent also occurs in deep neural networks \citep{belkin2019reconciling, spigler2019jamming, nakkiran2021deep} and random forests \citep{belkin2019reconciling}.  
 
Double descent directly challenges conventional statistical thinking  by suggesting that fitting large models in the overparameterized regime ($d >n$) is better than trying to fit optimal sparse models in the underparameterized regime ($d < n$). 
Confusingly, double descent also occurs with increasing $n$ \citep[e.g.][]{nakkiran2021deep}.  As this represents moving from the overparameterized to the underparameterized regime, it suggests the opposite conclusion, namely that for large datasets, fitting small models is better than fitting large models.

\subsection{Our setting}

This study aims to improve our understanding of sparsity and overfitting when fitting linear models by exploring both empirically and theoretically the effect of including either irrelevant predictors (increasing $d$) or irrelevant observations (increasing $n$). Specifically,  suppose that we arrange the training data in a response vector $\yy=[y_1,\ldots,y_{n_0}]^T \in \mathbb{R}^{n_0}$ and a $n_0\times d_0$ matrix of predictors $\mathbf{X}=[\xx_1,\ldots,\xx_{n_0}]^T$, such that $\xx_i \sim \mbox{independent} N_{d_0}(\zer_{d_0},\II_{d_0})$, where $\zer_{d_0}$ is the $d_0$-vector of zeros and $\II_{d_0}$ is the $d_0\times d_0$ identity matrix, and $\yy = \XX\be_0 + \ee$ with $\mathbf{e} \sim N_{n_0}(\zer_{n_0}, \sigma^2\II_{n_0})$, where $\be_0\in \mathbb{R}^{d_0}$ is an unknown regression parameter and $\sigma^2 >0$ is an unknown variance parameter.  We consider two sequences of training sets $\mathcal{D}$: I (Adding predictors; $d$ increases with $n=n_0$ fixed) for $d\le d_0$, we set $\mathcal{D}^{(d)}=(\yy,\XX^{(d)})$, where $\XX^{(d)}$ contains the first $d$ columns of $\XX$, and for $d>d_0$, we set $\mathcal{D}^{(d)}=(\yy,[\XX, \ZZ])$, where $\ZZ$ is an $n_0 \times (d-d_0)$ matrix with independent $N_{d-d_0}(\zer_{d-d_0},\II_{d-d_0})$ rows that are independent of $\XX$ and $\ee$; II (Adding observations; $n$ increases with $d=d_0$ fixed) for $n \le n_0$, we set $\mathcal{D}_{(n)}= (\yy_{(n)},\XX_{(n)})$, where $\yy_{(n)}$ contains the first $n$ elements of $\yy$ and $\XX_{(n)}$ contains the first $n$ rows of $\XX$, and for $n>n_0$, we set $\mathcal{D}_{(n)}=([\yy^T, \zer_{n-n_0}^T]^T, [\XX^T, \WW^T]^T)$,  where $\WW$ is an $(n-n_0) \times d_0$ matrix with independent $N_{d_0}(\zer_{d_0},\II_{d_0})$ rows that are independent of $\XX$ and $\ee$. We follow \cite{kobak2020optimal} in setting the responses corresponding to $\WW$ equal to zero (the marginal mean response); a more complicated alternative would be to set them equal to $N_{n-n_0}(\zer_{n-n_0}, \sigma^2\II_{n-n_0})$ noise that is independent of $\ee$, $\XX$ and $\WW$.

As we are interested in the conflicting consequences of the two types of sequences, we study them separately rather than together (i.e. both $d, n \rightarrow \infty$) as in \cite{bartlett2020benign} and \cite{hastie2022surprises}. These authors make all the predictors important ($d=d_0$), although \cite{hastie2022surprises} introduced noise by adding measurement error to the predictors. \cite{nakkiran2021deep} increased the number of predictors and the number of observations separately, but made these all important predictors and true observations.  Intuitively, insofar as many predictors are ``less important'' (have non-zero but small coefficients), they have a similar effect to noise predictors.  However, clearly separating important and noise predictors as we do is essential for understanding their separate effects.  \cite{mitra2019understanding}, albeit studying different estimators from us, included noise predictors but mixed them with important predictors and masked some of the effects we observe.  \cite{kobak2020optimal} also considered both sequences I and II but used specially normalized noise, replacing $\ZZ$ by $(d-d_0)^{-1/2}\lambda_z^{1/2}\ZZ$ and $\WW$ by $(n-n_0)^{-1/2}\lambda_w^{1/2}\WW$, where $\lambda_z, \lambda_w > 0$.  We argue that our results for unnormalized noise predictors are more interesting and relevant because normalized noise predictors have to be constructed and added by the analyst whereas we allow noise predictors to simply occur in the data and/or to be added by the analyst. The noise predictors of \cite{kobak2020optimal} are also asymptotically zero, so increasingly negligible, and following the general recommendation to standardize the predictors (to have variance one) to reduce heterogeneity would remove the normalization.  \cite{muthukumar2020harmless} considered sequences similar to our sequence I, but in their empirical work they only considered one of the cases we consider, namely $d_0 << n$. In this case, we need to include a large number of noise predictors to reach the interpolation point $d=n$, a circumstance that hides some findings we make by considering a wider set of relationships between $d_0$ and $n$.  Finally, \cite{hellkvist2023estimation} also considered our sequence I (calling noise predictors fake features) when the unknown parameters are zero-mean random vectors.  They did not consider noise observations.

\subsection{Our results}

In our empirical and theoretical study of the estimators and their test errors, we establish two important sets of results.

1) Adding noise predictors or observations to the data as in sequences I and II shrinks the estimators to zero as $d\rightarrow \infty$ or $n\rightarrow \infty$.   We compute the test errors for the estimators under the two sequences and prove that they both asymptote to the same value as $d\rightarrow 0$ or $\infty$ and $n\rightarrow 0$ or $\infty$. This is a key result for understanding double descent: the convergence of test errors to the same value at both zero and infinity, combined with the `condition number anomaly' (Dax 2022), which is reflected in the test error becoming infinite in Theorems \ref{thrm_seq1} and \ref{thrm_seq2}, collectively drive the double descent phenomenon.

Although there are results relevant to sequence I in the literature (e.g. similar expressions for test errors in \cite{hellkvist2023estimation}), sequence II has not been considered and our results for sequence II are new.  Even for sequence I, our results are not as obvious as they may appear.  For example, the estimators do not shrink to zero along the sequences of \cite{kobak2020optimal} or in the setting without noise predictors of \cite{belkin2020two} who suggest that using as many predictors as possible may be optimal in their setting. For a particular case of sequence I and with a different focus on alias and spuriously correlated predictors, \cite{muthukumar2020harmless}  discussed `\textit{signal bleeding}' into a large number of alias predictors -- which is over-shrinkage --  and \textit{overfitting of noise} under parsimonious selection of noise predictors in scenarios where sparse or limited data makes it challenging to distinguish between true predictors and those that are spuriously correlated with the outcome in the training set.   For the test errors when increasing the number of predictors, both \cite{muthukumar2020harmless} and \cite{hastie2022surprises} found as we do that the test error of the model approaches that of the null model.  We show empirically and theoretically that this holds in greater generality than considered by \cite{muthukumar2020harmless}, when we add noise predictors (rather than predictors with measurement error as in \cite{hastie2022surprises}), and for sequence II.   We suggest that in a limited training set, the less-important predictors in \cite{bartlett2020benign} behave similarly to noise, shrinking coefficient estimates toward zero. Importantly, this implies that overfitting noise is not always harmless or benign \citep{muthukumar2020harmless, bartlett2020benign, tsigler2023benign}.  Finally, for well-specified models, \cite{nakkiran2019more} related the peak in the test error to the condition number of the data matrix.  In addition, we use the ``condition number anomaly''  \citep{dax2022smallest} to obtain additional theoretical insights into why the peak always occurs.  This result combined with the fixed asymptotes of the test error establishes that double descent has to occur in our setup (and incidentally, also in that of \cite{hastie2022surprises}). Showing this for both sequence I and II is important as it highlights the importance of shrinkage induced by noise rather than the specific form of the noise in driving double descent.


2) For both sequences I and II, the estimators are approximately ridge regression estimators.  Moreover, for sequence II, the ridge regression estimator applied to the augmented data is also approximately a ``double shrunk'' ridge regression estimator.  The ``double shrinking'' means that the optimal value of the ridge parameter in this estimator can be negative when $n \rightarrow \infty$, even when the predictors are independent and the number of predictors $d=d_0$ is small (i.e. a low-dimensional setting).

The interpretation of least squares estimators applied to noisy data as being like ridge regression estimators applied to ``pure'' data originated in the measurement error setting \citep{webb1994functional, bishop1995training, hastie2022surprises}.   This is different from adding additional noise predictors or observations.  \cite{kobak2020optimal} added additional normalized noise predictors, but ignored the estimated coefficients of the noise predictors, and noted that a similar interpretation holds when we add normalized noise observations to the data. Our results for the full regression parameter estimator avoid the need to normalize the noise predictors and observations differently from the true versions.  The seminal paper on ridge regression \citep{hoerl1970ridge} showed that for $d<n$, a ridge estimator with a positive ridge parameter $\lambda$ has a lower mean squared error than the OLS estimator, provided the predictor matrix is of full rank $d$. Subsequent research  by \cite{dobriban2018high} and \cite{hastie2022surprises} in the high dimensional setting showed that, for predictors with any covariance matrix $\boldsymbol{\Sigma}$, as $d, n \to \infty$ with $d/n=\gamma$, the optimal ridge parameter $\lambda_{\text{opt}}$ is positive when the orientation of $\boldsymbol{\beta}$ is random. \cite{nakkiran2020optimal} showed that a positive $\lambda$ is optimal regardless of the rank of the predictor matrix or the orientation of $\boldsymbol{\beta}$, as long as the predictors are independent.  Recently,  \cite{liang2020just} and \cite{kobak2020optimal} have shown that $\lambda_{opt}$ can be zero or negative when the predictors are correlated. Furthermore, a recent non-asymptotic analysis by \cite{tsigler2023benign} showed that under a specific `spiked' covariance structure for the data in which $\boldsymbol{\Sigma}$ has a few large eigenvalues, a negative $\lambda$ can improve generalization bounds.  We find that $\lambda_{opt}$ is positive under sequence I (which is consistent with \cite{nakkiran2020optimal}) but can be negative under sequence II, even when the predictors are independent and even in  the low dimensional setting. This is a completely new result.


\bigskip
Our results under sequence I show double descent occurring in the overparameterized regime, but our results for sequence II show it occurring in the underparameterized regime. This shows that double descent is not driven directly by overparameterization, but rather is driven by shrinkage which can be induced by overparameterization due to including irrelevant noise predictors or by including irrelevant noise observations (or by measurement error in the predictors). The emphasis should be placed on choosing the correct shrinkage, rather than on fitting the exactly correct (possibly sparse) model ($d=d_0$) or an overparameterized (usually complex) model ($d>>d_0$).

The remainder of the paper is organised as follows.  We present numerical results to illustrate our key findings in Section \ref{sec:results}  and then present theoretical results that explain our numerical results in Section \ref{sec:theory}.  We present a data application in Section \ref{application} before concluding with a discussion in Section \ref{sec:conc}. Additional material and proofs of our theorems are given in Supplementary Material.


\section{Empirical results}
\label{sec:results}

In this section, we empirically explore the effect of noise in the form of noise predictors or noise observations, and its relationship to double descent.  We intentionally consider linear regression to make the set up as simple as possible, examine the effect of increasing $d$ or $n$, and then isolate and interpret the observed effects.

Our first set of experiments explores the effect of adding an increasing number of irrelevant noise predictors (i.e. increasing $d$) for different fixed values of $d_0$ with $n=n_0=50$.  Full details of the data generating steps and constructing the figures are set out in Algorithm S1 in the Supplementary Material.

Figure \ref{DDcurvesWS} shows the effect on the test error of increasing $d$ by adding noise predictors.  We observe a second descent in all cases, but the second descent does not always correspond to the optimal model in test error. In Figure \ref{DDcurvesWS}a and \ref{DDcurvesWS}d, the data are generated by a model with $d_0<n$ and an under-parameterized model performs better than the overparameterized model defined by the second descent. When $d_0$ is near the interpolation point ($d_0=n$), the optimal minimum test error is likely to occur during the second descent after escaping the zone around the interpolation point (see Figure \ref{DDcurvesWS}b and \ref{DDcurvesWS}e). If $d_0>n$, we obtain a useful second descent in the overparameterized regime in the sense that the global minimum of the test error is achieved in the second descent (see Figure \ref{DDcurvesWS}c and \ref{DDcurvesWS}f).  Clearly, we need to consider the dimension $d_0$ of the data generating model (also called the true dimension) as well as $d$ and $n$ to understand double-descent. It is also known that the signal-to-noise ratio of the model affects double descent \citep{hastie2022surprises}. We demonstrate this by recreating Figures \ref{DDcurvesWS}a, \ref{DDcurvesWS}b, and \ref{DDcurvesWS}c  with a higher signal-to-noise ratio in Figures \ref{DDcurvesWS}d, \ref{DDcurvesWS}e, and \ref{DDcurvesWS}f.  

\begin{figure}[t!]
\begin{center}$
\begin{array}{c}
\includegraphics[width=5cm, height=6cm, keepaspectratio]{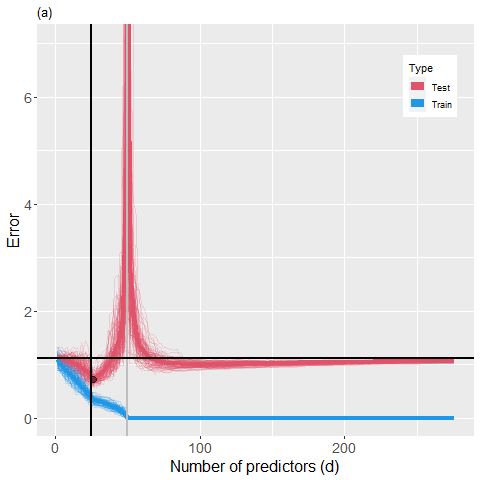}
\includegraphics[width=5cm, height=6cm, keepaspectratio]{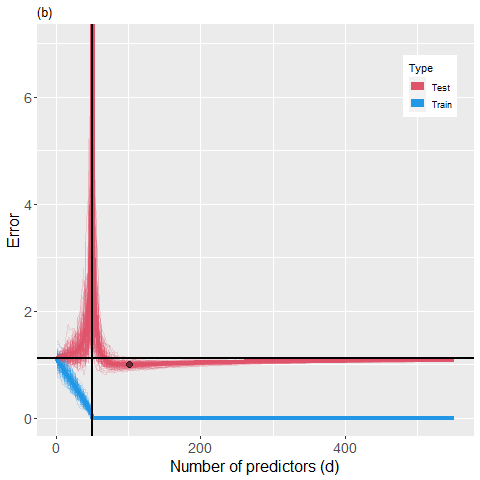}
\includegraphics[width=5cm, height=6cm, keepaspectratio]{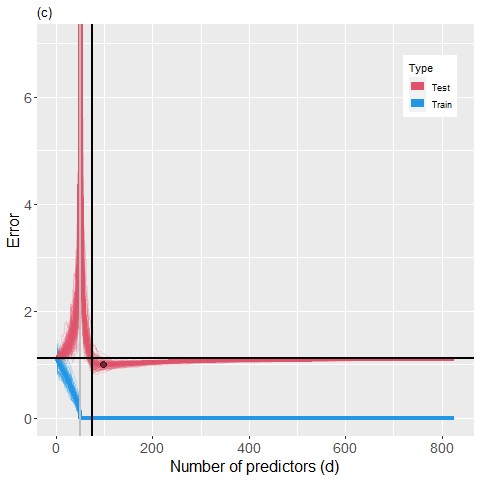}\\
\includegraphics[width=5cm, height=6cm, keepaspectratio]{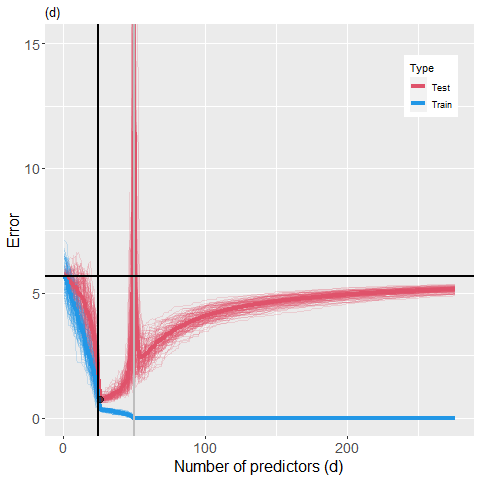}
\includegraphics[width=5cm, height=6cm, keepaspectratio]{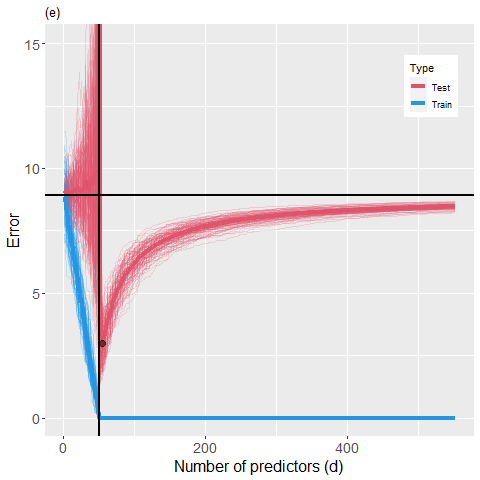}
\includegraphics[width=5cm, height=6cm, keepaspectratio]{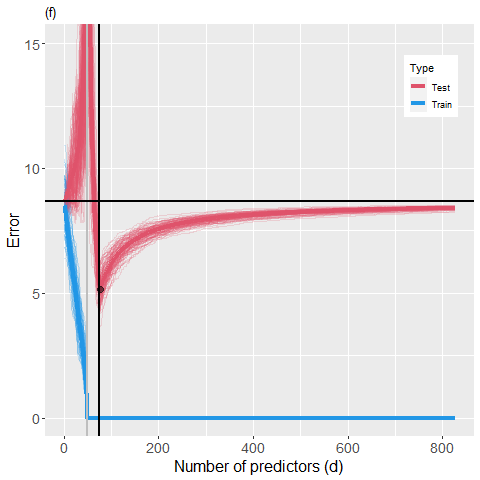}
\end{array}$
  \caption{Root-mean-squared training (blue) and test (red) errors obtained using Algorithm S1. We used $n_0=50$ and varied $d_0$ in the subplots as (a) and (d) $d_0=25$, (b) and (e) $d_0=50$, and (c) and (f) $d_0=75$. In subplots (a), (b), and (c), we have a weak signal-to-noise ratio ($\boldsymbol{\beta}_0=\boldsymbol{\beta}/\|\boldsymbol{\beta}\|$), while in subplots (d), (e), and (f), we have a strong signal-to-noise ratio ($\boldsymbol{\beta}_0=\boldsymbol{\beta}$). The vertical grey and black lines mark $n_0$ and $d_0$, respectively, and the horizontal black line marks the error of the null model with $\hat{\boldsymbol{\beta}}^{(0)}=\mathbf{0}$.} \label{DDcurvesWS}
\end{center}
\end{figure}

\begin{figure}[t!]
\begin{center}$
\begin{array}{c}
\includegraphics[width=5cm, height=6cm, keepaspectratio]{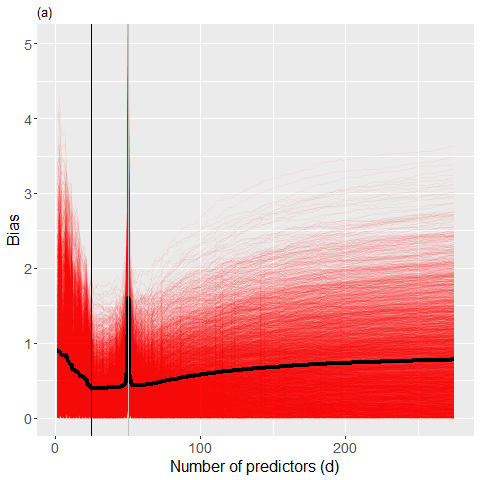}
\includegraphics[width=5cm, height=6cm, keepaspectratio]{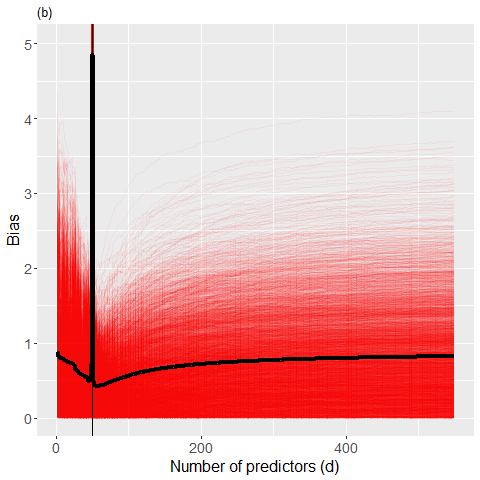}
\includegraphics[width=5cm, height=6cm, keepaspectratio]{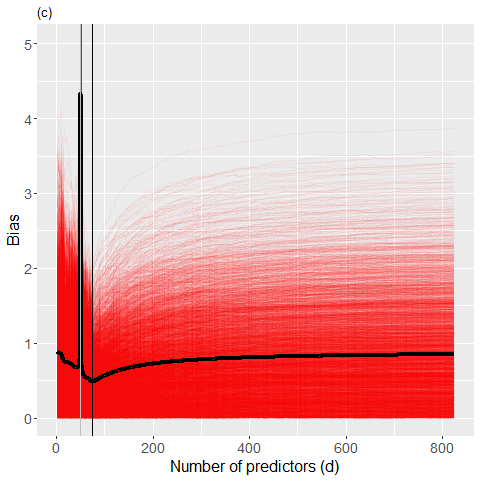}\\
\includegraphics[width=5cm, height=6cm, keepaspectratio]{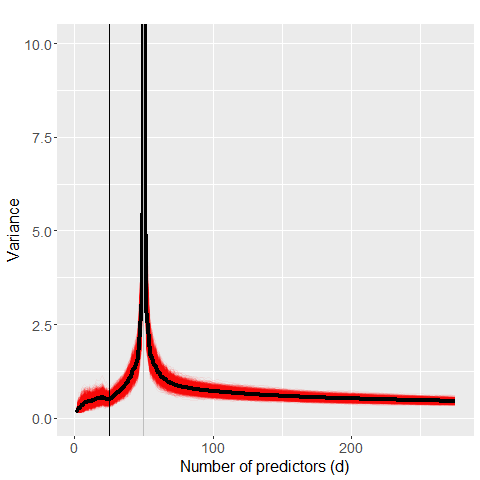}
\includegraphics[width=5cm, height=6cm, keepaspectratio]{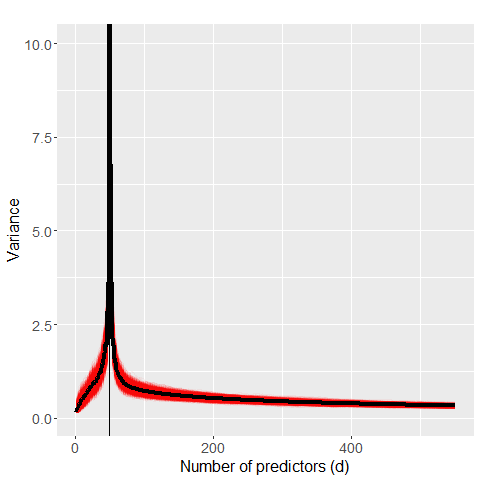}
\includegraphics[width=5cm, height=6cm, keepaspectratio]{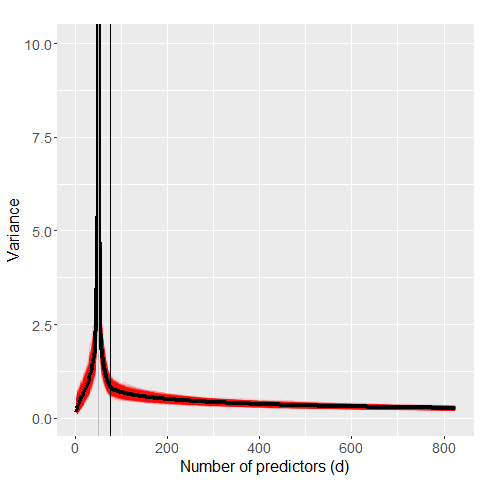}
\end{array}$
  \caption{Bias and variance components of the test error shown in Figure \ref{DDcurvesWS}. Row 1 shows the bias and row 2 shows the variance; columns (a), (b), and (c) correspond to the settings in subplots (a), (b), and (c) in Figure \ref{DDcurvesWS}, respectively. The individual red lines are the averages of 100 predictions for each of the 5000 test samples, and the bold black line is the average of 5000 red lines.}
  \label{BiasVarDec}
\end{center}
\end{figure}

The second descent in the overparameterized regime has led to questioning the classical bias-variance trade-off principle because bias and variance are considered irrelevant in the overparameterized regime. However, the test error can be decomposed into bias and variance terms in either regime, and we can still examine these functions.  The bias and variance components of the test error shown in Figures \ref{DDcurvesWS}a, \ref{DDcurvesWS}b, and \ref{DDcurvesWS}c are plotted in Figure \ref{BiasVarDec}. Overall, the bias starts out relatively high but decreases as $d$ increases. The bias reaches a minimum in the under-parameterized regime and increases to its highest value near the interpolation point. It reaches another minimum in the overparameterized regime after escaping the zone around the interpolation threshold and increases again, ultimately asymptoting to a limit. The variance, on the other hand, starts relatively low, but as the complexity of the model increases, it increases to its highest point near the interpolation point. Subsequently, in the overparameterized regime, the variance decreases, eventually reaching a minimum equal to the error variance ($\sigma^2=0.25$ for these data) asymptotically.  

\begin{figure}[t!]
\begin{center}$
\begin{array}{c}
\includegraphics[width=10cm, height=10cm, keepaspectratio]{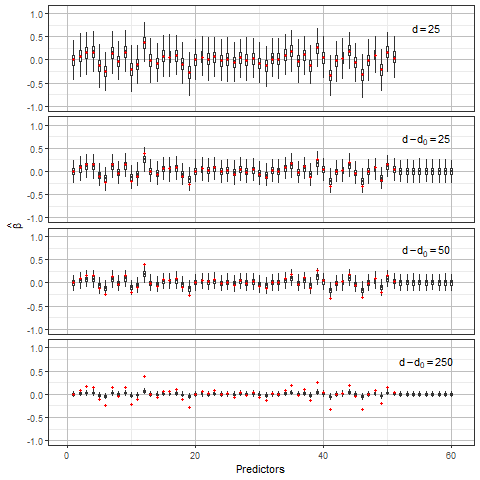}
\end{array}$
  \caption{Estimated regression coefficients from 1000 trials, with $n_0=50$ and $d_0=50$. These four subplots correspond to four values of $d$ along the horizontal axis in Figure \ref{DDcurvesWS}b. The red dots indicate the true nonzero coefficients; the true zero coefficients for the noise variables are not shown.  The estimated coefficients associated with the $d-d_0$ noise variables behave similarly, so we show the results only for the first nine noise variables. No noise variables ($d < d_0$) are in the top subplot, and 25 out of $d_0=50$ variables were randomly selected to predict the response. In the bottom three subplots, all $50$ variables were used together with $d-d_0 \in\{25, 50, 250\}$ noise variables. The first boxplot in each subplot represents the estimated intercept.} \label{BoxPlts}
\end{center}
\end{figure}

Examining the effect of including noise predictors on the parameter estimates in the working model is also interesting. Figure \ref{BoxPlts} shows boxplots of the estimated intercept and the slope coefficients up to the first 60 predictors included in each working model with $d_0=50$ and $n=n_0=50$ from $1000$ trials. In the top panel of Figure \ref{BoxPlts}, 25 important predictors are omitted. Thus, the coefficient estimates for the included predictors are unbiased (the medians coincide with the true values, as indicated by red dots). This occurs because the omitted predictors are independent of the important predictors included in the model; hence, there is no omitted variable bias. However, it is important to note that the test error includes bias due to the underfitted model, as shown in Figure \ref{BiasVarDec}a. In the bottom three panels of Figure \ref{BoxPlts}, the estimates become increasingly biased as we include more irrelevant noise predictors -- the estimated coefficients shrink towards zero. Larger coefficients shrink at a faster rate than smaller ones. The variances of the estimates increase as we approach the interpolation point in the underparameterized regime, peaking at the interpolation point (results not shown), and decrease as we add more irrelevant predictors in the overparameterized regime. Thus, the noise predictors have a shrinkage effect on the estimated coefficients, resulting in them approaching zero as $d-d_0$ increases. 

So far, we have treated the size of the training set $n$ as fixed and varied the number of predictors $d$ (corresponding to parameters) in the working model. It is also plausible that if $n$ is large, some of the observations in the training set may be noise observations that are not generated by the data generating model. Therefore, we consider the problem in which $d$ is fixed and $n$ increases as we add irrelevant noise observations to the training data.  In the additional observations, we set the response equal to zero for simplicity. Full details of the data generating steps and constructing the figures are set out in Algorithm S2 in the Supplementary Material.

Figure \ref{DDcurvesNObs} shows the effect on the training and test errors, and Figure \ref{BoxPltsRowNoise} shows boxplots of the regression coefficients as we increase $n$ from the overparameterized to the underparameterized regime.   In Figures \ref{DDcurvesNObs}(a), \ref{DDcurvesNObs}(d), and \ref{DDcurvesNObs}(e), the number of standard normal noise observations required to escape the zone near the interpolation point $n=d$ causes overshrinkage, as shown in Figure \ref{BoxPltsRowNoise}. Consequently, the minimum test error is achieved in the overparameterized regime.  Figure \ref{DDcurvesNObs}(b) illustrates that, by incorporating random observations, a global minimum in the underparameterized regime can be achieved with an appropriate amount of noise. Figures \ref{DDcurvesNObs}(c) and \ref{DDcurvesNObs}(f) show that when the signal-to-noise ratio is low, more noise is required to reach the lowest test error than when the signal-to-noise ratio is high. Interestingly, Figure \ref{DDcurvesNObs} reveals that for double descent, the number of predictors $d$ and the training sample size $n$ do not need to be large. Since double descent can occur for small $d_0$, it is not the overparameterization of the model with noise predictors that leads to improved prediction performance but the regularization due to the inclusion of noise that improves the model performance. 

\begin{figure}[t!]
\begin{center}$
\begin{array}{c}
\includegraphics[width=5cm, height=6cm, keepaspectratio]{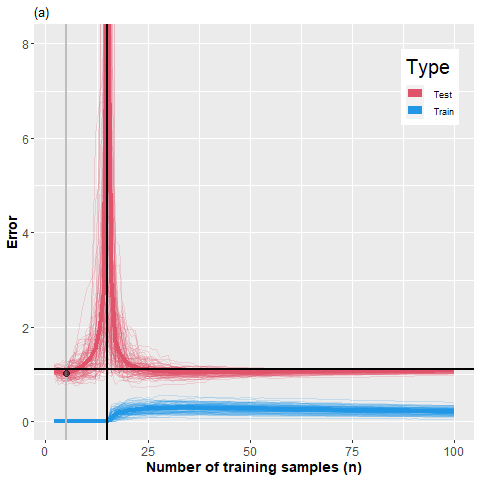}
\includegraphics[width=5cm, height=6cm, keepaspectratio]{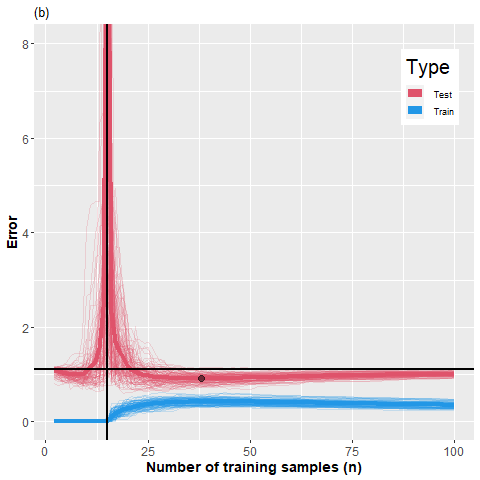}
\includegraphics[width=5cm, height=6cm, keepaspectratio]{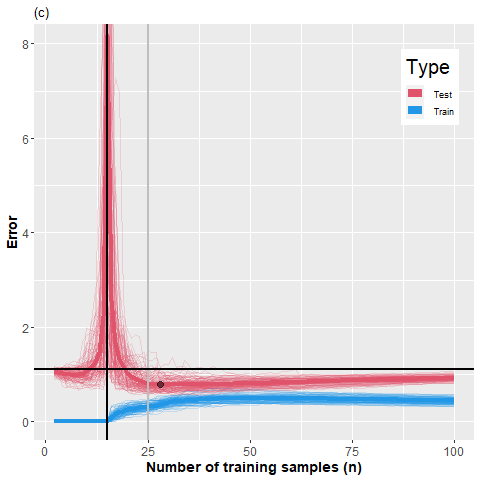}\\
\includegraphics[width=5cm, height=6cm, keepaspectratio]{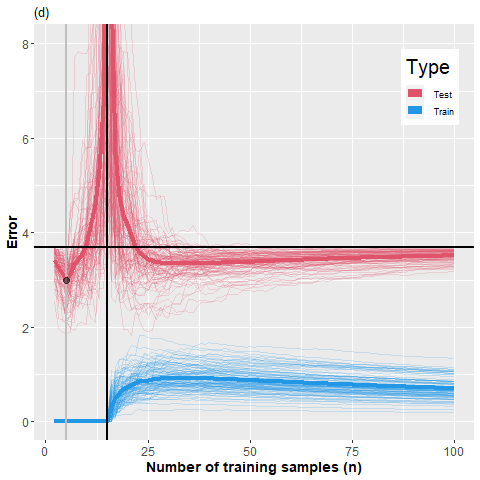}
\includegraphics[width=5cm, height=6cm, keepaspectratio]{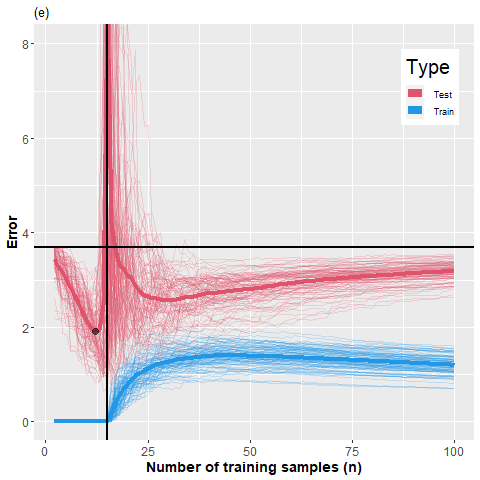}
\includegraphics[width=5cm, height=6cm, keepaspectratio]{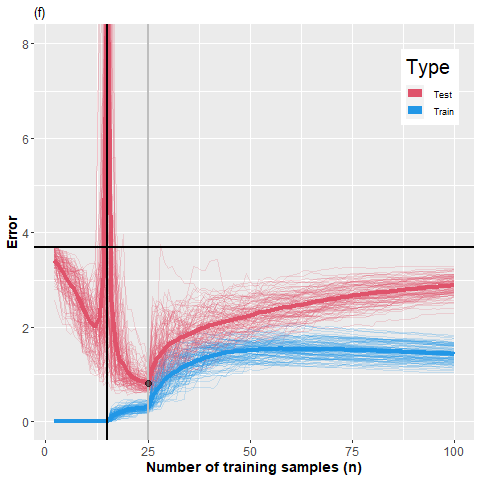}
\end{array}$
  \caption{Root-mean-squared errors obtained using Algorithm S2. We used $d=d_0=15$ and varied $n_0$, the number of observations from the model in the training, in subplots as (a) and (d) $n_0=5$, (b) and (e) $n_0=15$, and (c) and (f) $n_0=25$. In subplots (a), (b), and (c), we have a weak signal-to-noise ratio ($\boldsymbol{\beta}_0=\boldsymbol{\beta}/\|\boldsymbol{\beta}\|$), while in subplots (d), (e), and (f), we have a strong signal-to-noise ratio ($\boldsymbol{\beta}_0=\boldsymbol{\beta}$). The vertical grey and black lines mark $n_0$ and $d_0$, respectively, and the horizontal black line marks the error of the null model with $\hat{\boldsymbol{\beta}}_{(n)}=\mathbf{0}$.
}\label{DDcurvesNObs}
\end{center}
\end{figure}

\begin{figure}[t!]
\begin{center}$
\begin{array}{c}
\includegraphics[width=8cm, height=10cm, keepaspectratio]{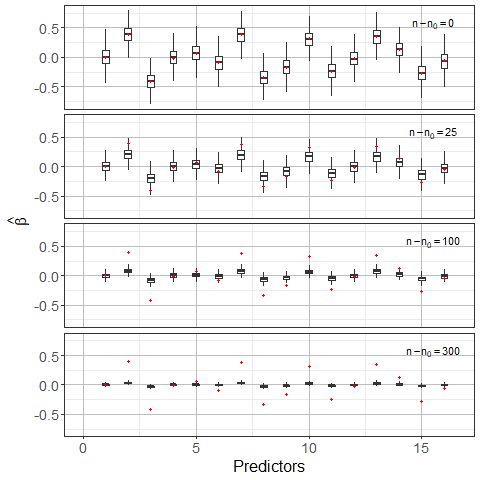}
\end{array}$
  \caption{Estimated regression coefficients from 1000 trials with $d_0=15$ and $n_0=25$ (the setting in Figure \ref{DDcurvesNObs}c). The four subplots correspond to four values of $n$ along the horizontal axis in Figure \ref{DDcurvesNObs}c.  The red dots indicate the true non-zero coefficients. The top subplot has no noise observations ($n=n_0$). The first boxplot in each subplot represents the estimated intercept. } \label{BoxPltsRowNoise}
\end{center}
\end{figure}

\begin{figure}[t!]
\begin{center}$
\begin{array}{c}
\includegraphics[width=8cm, height=10cm, keepaspectratio]{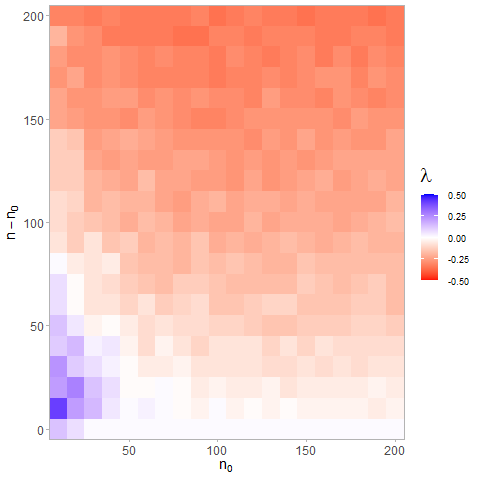}
\end{array}$
  \caption{Heatmap for the average optimal ridge parameter $\lambda$ over 500 trials as a function of $n_0$ and $n-n_0$. Here, ``optimal'' means minimizing the test error.  The data were generated following Algorithm S2 with $d_0=30$.  We set $\boldsymbol{\beta}_0=\boldsymbol{\beta}$ for a high signal-to-noise ratio. Purple-blue values are positive and orange-red values are negative.} \label{NegLambdaIdentityCov}
\end{center}
\end{figure}

We find that OLS/minimum norm OLS estimation performs regularization in the presence of noise. However, the amount of noise present in real data may not provide optimal regularization, resulting in undershrinkage or overshrinkage. Undershrinkage could be addressed by generating and adding artificial noise, but it seems easier to control the shrinkage by using explicit regularization as in ridge regression. If there is overshrinkage through too much noise in the model, then ridge regression with $\lambda > 0$ increases the shrinkage and makes prediction worse.  On the other hand, ridge regression with a negative $\lambda$ can reduce the shrinkage due to noise and, as shown by \cite{kobak2020optimal} in a particular scenario, the optimal $\lambda$ can be zero or negative. We contend that this phenomenon is not limited to such specific cases and can manifest in many situations when noise is present in the data.  Figure \ref{NegLambdaIdentityCov} illustrates that the optimal value of $\lambda$ can be negative when there are large numbers of noise observations.

It has been established numerically in the literature that more than two descents can occur in particular settings \citep{nakkiran2020optimal,liang2020multiple}. For independent predictors, we found that multiple descents occur when the predictors are not on the same scale and that multiple descents can be reduced by ensuring that all predictors are on the same scale; however, double descent remains inevitable. This is shown numerically in Figure S1 in the Supplementary Material. 



\section{Theoretical results}
\label{sec:theory}


\subsection{The estimators}
\label{subsec:estimators}
We construct the two sequences of training data $\mathcal{D}^{(d)}$ and $\mathcal{D}_{(n)}$ as described in the Introduction.  Then we fit working linear regression models to the two sequences of training data by estimating the regression parameters using the ordinary least squares (OLS) and minimum norm OLS estimators when $d < n_0$ and $d > n_0$ for sequence I, and $n>d_0$ and $n<d_0$ for sequence II, respectively.  We write these linear estimators as $\hat{\be}^{(d)} =\MM^{(d)}\yy$ for sequence I and $\hat{\be}_{(n)} =\MM_{(n)}\yy$ for sequence II, where $\MM^{(d)}$ is a $d\times n_0$ matrix and $\MM_{(n)}$ is a $d_0\times n_0$ matrix. The matrix $\MM^{(d)}$ depends on $\XX^{(d)}$ or $[\XX,\ZZ]$ and $\MM_{(n)}$ depends on $\XX_{(n)}$ or $[\XX^T,\WW^T]^T$, but neither depends on $\yy$. For sequence I, when $d \le\min(n_0,d_0)$,  $\MM^{(d)} = (\XX^{(d)T}\XX^{(d)})^{-1}\XX^{(d)T}$.  If $d_0 >n_0$, we switch for $n_0 < d \le d_0$ to the minimum norm OLS estimator with $\MM^{(d)} = \XX^{(d)T}(\XX^{(d)}\XX^{(d)T})^{-1}$ and then for $d > d_0$  $\MM^{(d)} = [\XX, \ZZ]^{T}(\XX\XX^{T} + \ZZ\ZZ^T)^{-1}$; if $d_0 < n_0$, for $d_0 < d \le n_0$, we use $\MM^{(d)} = \{[\mathbf{X}, \mathbf{Z}]^T[\mathbf{X}, \mathbf{Z}]\}^{-1} [\mathbf{X}, \mathbf{Z}]^T$
before we switch to the minimum norm OLS estimator for $d \ge n_0$.  For sequence II, when $n \le \min(n_0,d_0)$,  $\MM_{(n)} = \XX_{(n)}^T(\XX_{(n)}\XX_{(n)}^T)^{-1}[\II_n,\zer_{n,n_0-n}]$.  If $d_0 < n_0$, we switch to the OLS estimator for $d_0 < n \le n_0$ so $\MM_{(n)} = (\XX_{(n)}^T\XX_{(n)})^{-1}\XX_{(n)}^T[\II_n,\zer_{n,n_0-n}]$ and then for $n>n_0$ $\MM_{(n)} = (\XX^T\XX +\WW^T\WW)^{-1}\XX^T$; if $d_0 > n_0$, for $ n_0 < n \le d_0$ $\MM_{(n)} = \XX^T\bA^{11} + \WW^T\bA^{21}$, where $\bA^{11} = \{\XX\XX^T - \XX\WW^T(\WW\WW^T)^{-1}\WW\XX^T\}^{-1}$ and $\bA^{21} = -(\WW\WW^T)^{-1}\WW\XX^T\bA^{11}$, before we switch to the OLS estimator for $n>d_0$. The second expression for $\MM_{(n)}$ is obtained by letting $\bA = [\XX^T,\WW^T]^T[\XX^T,\WW^T]$, partitioning it into its natural four submatrices and then using standard formulas for the inverse of a partitioned matrix to obtain the corresponding terms in $\bA^{-1}$.

We showed in Figures \ref{BoxPlts} and \ref{BoxPltsRowNoise} that adding noise predictors in sequence I or noise observations in sequence II induces shrinkage in $\hat{\be}^{(d)}$ and $\hat{\be}_{(n)}$, respectively, and eventually shrinks these estimators to zero.  We now confirm this effect theoretically as $d$ or $n \rightarrow \infty$.  We let `a.s.' mean `almost surely' and `P' mean 'in probability'.  Also, let $\|\xx\|=(\xx^T\xx)^{1/2}$ denote the Euclidean norm of $\xx$.  For sequence I, for $d > \max(d_0,n_0)$, write $\hat{\boldsymbol{\beta}}^{(d)} = (\hat{\boldsymbol{\beta}}^{(d)T}_{d_0}, \hat{\boldsymbol{\beta}}^{(d)T}_{-d_0})^T$, where $\hat{\boldsymbol{\beta}}^{(d)}_{d_0}= \mathbf{X}^T(\mathbf{X}\mathbf{X}^T+\ZZ\ZZ^T)^{-1}\mathbf{y}$ denotes the first $d_0$ elements and $\hat{\boldsymbol{\beta}}^{(d)}_{-d_0}=\mathbf{Z}^T(\mathbf{X}\mathbf{X}^T+\ZZ\ZZ^T)^{-1}\mathbf{y}$ denotes the remaining $d-d_0$ elements of $\hat{\boldsymbol{\beta}}^{(d)}$.   By the strong law of large numbers, the $n_0\times n_0$ matrix $\XX\XX^T + \ZZ\ZZ^T  = O(d)$ a.s. as $d \rightarrow \infty$.  Since $\hat{\boldsymbol{\beta}}^{(d)}_{d_0}\in \mathbb{R}^{d_0}$ and $\hat{\boldsymbol{\beta}}^{(d)}_{-d_0}\in \mathbb{R}^{d-d_0}$, we have  $\|\hat{\boldsymbol{\beta}}^{(d)}_{d_0}\| = O(d^{-1}d_0^{1/2})$ and $\|\hat{\boldsymbol{\beta}}^{(d)}_{-d_0}\| =O(d^{-1}(d-d_0)^{1/2})$ a.s. as $d \rightarrow \infty$.  Similarly, for sequence II, for $n > \max(n_0,d_0)$, $\hat{\boldsymbol{\beta}}_{(n)}= (\mathbf{X}^T\mathbf{X}+\WW^T\WW)^{-1}\mathbf{X}^T\mathbf{y}$, the $d_0\times d_0$ matrix $\XX^T\XX+\WW^T\WW=O(n)$ a.s. as $n \rightarrow \infty$, and  $\|\hat{\be}_{(n)}\|=O(n^{-1}n_0)$ a.s. as $n \rightarrow \infty$. That is, for $d_0$ and $n_0$ fixed or increasing sufficiently slowly in sequence I and II, respectively, noise predictors or observations eventually shrink the estimated coefficients to zero.

\subsection{The test error}

We examine the performance of the working models fitted along the two training sequences through the test error.  In our empirical work, we compute empirical test errors by computing the mean squared error over a large test sample.  In theoretical work, we calculate the model-expectation of the squared error at a single observation $(y', \xx^{'(d)}) \in \mathbb{R}^{d+1}$ that is generated independently from the same data generating model as the training data, namely $\xx' \sim N_{d_0}(\zer_{d_0},\II_{d_0})$, $y' = \xx^{'T}\boldsymbol{\beta}_0 + e'$ with $e' \sim N(0, \sigma^2)$, and $\xx^{'(d)}$ contains the first $d$ elements of $\xx'$ if $d \le d_0$, and $\xx^{'(d)}=[\xx^{'T}, \zz^{'T}]$, where $\zz' \sim N_{d-d_0}(\zer_{d-d_0},\II_{d-d_0})$ independent of $\xx'$, if $d >d_0$.  We obtain different test errors, depending on what we condition on.  For example, corresponding to the light red curves in our figures, we have the conditional (on the training data) test errors
\[
R^{(d)}(\mathcal{D}^{(d)}) = \left\{ \begin{array}{ll}\mbox{E}_{y',\mathbf{x}'}(y' - \mathbf{x}^{'(d)T} \hat{\boldsymbol{\beta}}^{(d)})^2 & d\le d_0\\ \mbox{E}_{y',\mathbf{x}'}(y' - [\mathbf{x}^{'T},\zz^{'T}] \hat{\boldsymbol{\beta}}^{(d)})^2& d>d_0 \end{array} \right ., \qquad R_{(n)}(\mathcal{D}_{(n)}) = \mbox{E}_{y',\mathbf{x}'}(y' - \mathbf{x}^{'T} \hat{\boldsymbol{\beta}}_{(n)})^2.
\]
Note that for sequence II, we compute the test error under the model for the real data $(y', \mathbf{x}')$ rather than for noise observations $(0, \mathbf{w}')$ because this is the prediction we are actually interested in;  alternatives such as the test error under the noise observations or a linear combination of these two possibilities (e.g. with weights $n_0/n$ and $(n-n_0)/n$) could also be considered.
The dark red curves in our figures are (unconditional) test errors $R^{(d)}$  and $R_{(n)}$ obtained by taking the expectations (over $\mathcal{D}^{(d)}$ or $\mathcal{D}_{(n)}$) of $R^{(d)}(\mathcal{D}^{(d)})$ and $R_{(n)}(\mathcal{D}_{(n)})$.  

The test error under sequence I has been obtained by \cite{breiman1983how} and \cite{belkin2019reconciling}.  We state the theorem for completeness.
\begin{thrm} \citep{breiman1983how, belkin2019reconciling} \label{thrm_seq1}
Suppose that the training set is constructed as in sequence I.   If $d < d_0$, partition $\be_0$ into the first $d$ elements and the remaining $d_0-d$ elements as $\be_0 = [\be_{0}^{(d)T}, \tilde\be_{0}^{(d)T}]^T$; if $d \ge d_0$,  $\be_{0}^{(d)} = \be_0$ and $\tilde\be_{0}^{(d)}$ is replaced by zero.
Then, for $d \le n_0-2$, the estimator is OLS and the test error is
\begin{align*}
    R^{(d)}&= (\|\tilde\be_{0}^{(d)}\|^2 + \sigma^2)\Big(1 + \frac{d}{n_0-1-d}\Big), 
\end{align*}
for $d\in\{n_0-1,n_0,n_0+1\}$, the test error is infinite, and
for $d \ge n_0+2$, the estimator is minimum norm OLS and the test error is
\begin{align} \label{eq:thm4.1}
   R^{(d)} &= (\|\tilde\be_{0}^{(d)}\|^2 + \sigma^2)\Big(1 + \frac{n_0}{d-1-n_0}\Big) + \|\be_{0}^{(d)}\|^2\Big(1 - \frac{n_0}{d}\Big).
\end{align}

\end{thrm}
\noindent One insight that is not obvious from the form of Theorem \ref{thrm_seq1} is that for $d > \max(n_0,d_0)$, the test error $R^{(d)}$ in (\ref{eq:thm4.1}) can be written as the sum of $ \sigma^2 + \E\|\{\XX^T(\XX\XX^T+\ZZ\ZZ^T)^{-1}\XX -\II_{d_0}\}\be_0\|^2 + \sigma^2\E\trace\{\XX\XX^T(\XX\XX^T+\ZZ\ZZ^T)^{-2}\}$, the test error due to the first $d_0$ components of $\hat{\be}^{(d)}$, and $\E\|\ZZ^T(\XX\XX^T+\ZZ\ZZ^T)^{-1}\XX \be_0 \|^2 + \sigma^2\E\trace\{\ZZ(\XX\XX^T+\ZZ\ZZ^T)^{-2}\}$, the test error due to the last $d-d_0$ components of $\hat{\be}^{(d)}$, respectively.  See the Supplementary material for details.  We do not have analytic expressions for these component terms.  Nonetheless, they are important because they show that the estimated coefficients of the noise predictors contribute to the test error so these predictors need not be beneficial or even harmless.

The test error under sequence II is more complicated to work with than that under sequence I because analytic expressions are not available once we augment the data with noise observations ($n > n_0$).  Nonetheless, we obtain the following theorem.
\begin{thrm} \label{thrm_seq2}
Suppose the training set is constructed as in sequence II.  For $n <\min(n_0,d_0)$, the estimator is minimum norm OLS and the test error is
\begin{align*}
R_{(n)}&= \left\{\begin{array}{ll}\|\be_{0}\|^2\Big(1-\frac{n}{d_0}\Big) + \sigma^2\Big(1+\frac{n}{d_0-1-n}\Big),  & n \le d_0-2 \\\infty, & n \in \{d_0-1, d_0\}\end{array} \right. .
\end{align*}
If $n_0 < d_0$ and $n_0 \le n < d_0$, the estimator is minimum norm OLS so the test error is
\begin{align*}
R_{(n)}& = \sigma^2 + \E\|(\XX^T\bA^{11}\XX + \WW^T\bA^{21}\XX - \II_{d_0})\be_0\|^2\\
&+ \sigma^2\E\trace\{\bA^{11}\XX\XX^T\bA^{11} + (\bA^{21})^T\WW\XX^T\bA^{11} + \bA^{11}\XX\WW^T\bA^{21} + (\bA^{21})^T\WW\WW^T\bA^{21}\}.
\end{align*}
where $\bA^{11} = \{\XX\XX^T - \XX\WW^T(\WW\WW^T)^{-1}\WW\XX^T\}^{-1}$ and $\bA^{21} = -(\WW\WW^T)^{-1}\WW\XX^T\bA^{11}$.\\
If $d_0 < n_0$ and $d_0 \le n < n_0$, the estimator is OLS and the test error is
\begin{align*}
R_{(n)} &= \left\{\begin{array}{ll} \infty, & n \in \{d_0, d_0+1\} \\ \sigma^2\Big(1 +\frac{d_0}{n-d_0-1}\Big), & n \ge d_0+2\end{array}\right . .
\end{align*}
If $n > \max(n_0,d_0)$, the training set is augmented by noise observations, and the estimator is OLS, so the test error is
\begin{align} \label{eq:thm4.2}
R_{(n)} &= \sigma^2 + \E\|\{(\XX^T\XX + \WW^T\WW)^{-1}\XX^T\XX -\II_{d_0} \} \be_0\|^2  \nonumber\\
&\qquad + \sigma^2\E\trace\{\XX^T\XX(\XX^T\XX + \WW^T\WW)^{-2}\}.  
\end{align}
\end{thrm}
\noindent  The proof is given in the Supplementary material. As $d=d_0$, there is no decomposition of the test error (\ref{eq:thm4.2}) corresponding to that of (\ref{eq:thm4.1}).

\subsection{The effect of including noise on the test error}


We saw empirically in Figures \ref{BoxPlts} and  \ref{BoxPltsRowNoise}
and theoretically in Section \ref{subsec:estimators} that adding noise to the training data shrinks the estimators to zero as $d$ or $n \rightarrow \infty$. This means that, as shown in Figure \ref{DDcurvesWS}, the test errors should converge as $d$ or $n \rightarrow \infty$ to the test error of the null model with $d=n=0$. That is, both ends of the test error should equal the same value $\E(y^{'2}) = \|\be_0\|^2 + \sigma^2$.  This has also been observed by \cite{muthukumar2020harmless} and, in a different setup, by \cite{hastie2022surprises} but it has not been proved to occur.  That these asymptotes hold under sequence I follows immediately from Theorem \ref{thrm_seq1}.  It is much more challenging to prove that they hold under sequence II (because there is no analytic form for the test error under sequence II).  We now prove that this occurs in Theorem \ref{thrm_asmpt_test}; the proof is given in the Supplementary Material.

\begin{thrm}\label{thrm_asmpt_test}
Under sequence I, as $d \rightarrow \infty$ with $d_0$ fixed, 
$$
R^{(d)} = \sigma^2  + \|\be_0\|^2 +O(d^{-1}).
$$
Under sequence II, as $n\rightarrow \infty$ with $n_0$ fixed,
$$
R_{(n)} = \sigma^2 + \|\be_0\|^2+ O(n^{-1}).
$$
\end{thrm}
\noindent We can extend the proof under sequence II to obtain further terms in the approximation.  As these do not give any additional insight,  we have kept the result to the simplest one we need and use.



\subsection{Why double descent occurs}

The fact that both ends of the test error equal the same value means double descent occurs whenever the test error increases above the initial test error.  The test error can increase (e.g. Figure \ref{DDcurvesWS}c) or decrease (e.g. Figure \ref{DDcurvesWS}a) from its value at the origin. If the test error increases immediately, it has a local minimum at the origin (which we still call the first descent), while if it decreases and then increases above its initial value,  it has a local minimum before the interpolation point ($d=n$).  In either case, it will have to descend again (the second descent) to achieve its asymptote.  If it decreases below the null model test error, the test error has to increase again to achieve its asymptote.  
The second descent does not necessarily achieve a global minimum; whether this occurs or not is related to the amount of shrinkage induced by the noise variables or observations.
    
What makes the test error increase above the null model test error?  The estimators and their test errors involve $(\DD^T\DD)^{-1}$ or $(\DD\DD^T)^{-1}$ with $\DD$ equal to $\XX^{(d)}$, $[\XX,\ZZ]$, $\XX_{(n)}$ or $[\XX^T, \WW^T]^T$. The underlying driver behind the test error behavior is the generalized condition number $\kappa(\mathbf{D})$ of these matrices.  In our examples, the elements of the matrix $\mathbf{D}\in\mathbb{R}^{n\times d}$ are independently sampled from the same probability distribution. The generalized condition number $\kappa(\mathbf{D})$ of a $n\times d$ rectangular matrix $\mathbf{D}$ of rank $r \leq \mbox{min}(n,d)$ is $\kappa(\mathbf{D})=\gamma_1(\mathbf{D})/\gamma_r(\mathbf{D})$, where $\gamma_1(\mathbf{D}) \ge \ldots \ge \gamma_r(\mathbf{D})$ are the non-zero singular values of $\mathbf{D}$ in descending order \citep{zielke1988some}. As a direct consequence of the Cauchy Interlacing Theorem \citep[][Theorem 4.3.4]{horn2012matrix}, for a fixed $n$ and $d < n$, the condition number of $\mathbf{D}$ increases monotonically as $d\to n$; if $\mathbf{D}$ is of full rank, it reaches its maximum value when $d=n$. The increasing condition number of $\mathbf{D}$ causes the variance and the bias of the estimated coefficients to increase, with the maximum variance and maximum bias and hence maximum test error at the interpolation point.  If this maximum is above the asymptote, we observe a second descent in the overparameterized regime.  For $d > n$, the condition number of $\mathbf{D}$ is itself likely to decrease as we increase $d$ \citep[the``condition number anomaly'',][]{dax2022smallest}. This behavior is illustrated  in Figure S2 in the Supplementary Material.   The shape of the curves is strikingly similar to that of the variance curves in the second row of Figure \ref{BiasVarDec}.

\subsection{The inclusion of noise and ridge regression}

We can obtain additional insight into the estimators by approximating them by the ridge regression estimator $\hat{\boldsymbol{\beta}}_\lambda=(\mathbf{X}^T\mathbf{X}+\lambda \mathbf{I}_{d_0})^{-1}\mathbf{X}^T\mathbf{y}=\XX^T(\mathbf{X}\mathbf{X}^T+\lambda \mathbf{I}_{n_0})^{-1}\yy$. For sequence I, write $\hat{\boldsymbol{\beta}}^{(d)}_{d_0}-\hat{\boldsymbol{\beta}}_\lambda=\mathbf{X}^T\LL^{(d)}(\lambda)\mathbf{y}$, where
$\LL^{(d)}(\lambda) =(\mathbf{X}\mathbf{X}^T+\ZZ\ZZ^T)^{-1}(\lambda\mathbf{I}_{n_0}-\ZZ\ZZ^T)(\mathbf{X}\mathbf{X}^T+\lambda\mathbf{I}_{n_0})^{-1}$ is an $n_0\times n_0$ matrix.
Let $d^{-1}d_0 \rightarrow \delta$, with $0 \le \delta < 1$, as $d \rightarrow \infty$.  Then $(d-d_0)^{-1}\ZZ\ZZ^T -\II_{n_0} = O_p((d-d_0)^{-1/2})$ as $d \rightarrow \infty$ and hence, 
$\LL^{(d)}(d-d_0) = O_p(d^{-3/2})$, so $\|\hat{\boldsymbol{\beta}}^{(d)}_{d_0}-\hat{\boldsymbol{\beta}}_{d-d_0}\|=O_p(d^{-1})$  and, from above, $\|\hat{\boldsymbol{\beta}}^{(d)}_{-d_0}\|=O(d^{-1/2})$ a.s., so all the elements converge to zero. For sequence II,  let $n^{-1}n_0 \rightarrow\nu$ with $0 \le \nu < 1$. Then write $\hat{\boldsymbol{\beta}}_{(n)}-\hat{\boldsymbol{\beta}}_{(n_0),n-n_0}= \LL_{(n)}(\lambda)\XX^T\yy$, where the $d_0 \times d_0$ matrix $\LL_{(n)}(\lambda) = (\mathbf{X}^T\mathbf{X}+\WW^T\WW)^{-1}\{\lambda \mathbf{I}_{d_0}-\WW^T\WW\}\{\XX^T\XX + \lambda\II_{d_0}\}^{-1}$.
We have $\LL_{(n)}(n-n_0) = O_p(n^{-3/2})$ and $\|\hat{\boldsymbol{\beta}}_{(n)}-\hat{\boldsymbol{\beta}}_{n-n_0}\| = O_p(n^{-1/2})$. 
We gather the results into a Theorem.
\begin{thrm}\label{thrm_kobak}
Let $\hat{\boldsymbol{\beta}}_\lambda=(\mathbf{X}^T\mathbf{X}+\lambda \mathbf{I}_{d_0})^{-1}\mathbf{X}^T\mathbf{y}$ be the ridge regression estimator applied to the data $(\yy,\XX)$, let $\hat{\boldsymbol{\beta}}^{(d)}_{d_0}$ denote the first $d_0$-elements and $\hat{\boldsymbol{\beta}}^{(d)}_{-d_0}$ the remaining $d-d_0$ elements of $\hat{\boldsymbol{\beta}}^{(d)}$.  Under sequence I,  if $d^{-1}d_0 \rightarrow \delta$, with $0 \le \delta < 1$, as $d \rightarrow \infty$, then
 $\|\hat{\boldsymbol{\beta}}^{(d)}_{d_0}-\hat{\boldsymbol{\beta}}_{d-d_0}\|=O_p(d^{-1})$  and $\|\hat{\boldsymbol{\beta}}^{(d)}_{-d_0}\|=O_p(d^{-1/2})$.  Under sequence II,  if $n^{-1}n_0 \rightarrow\nu$ with $0 \le \nu < 1$, then $\|\hat{\boldsymbol{\beta}}_{(n)}-\hat{\boldsymbol{\beta}}_{(n_0), n-n_0}\| = O_p(n^{-1/2})$ as $n \rightarrow \infty$.
\end{thrm}

Related results were given by \cite{kobak2020optimal} for normalized noise variables.  Specifically, they normalized the noise predictors and replaced $\ZZ$ in the training set $\mathcal{D}^{(d)}$ by $\tilde{\ZZ}=(d-d_0)^{-1/2}\lambda_z^{1/2}\ZZ$, and replaced $\WW$ in the training set $\mathcal{D}_{(n)}$ by $\tilde{\WW}=(n-n_0)^{-1/2}\lambda_w^{1/2}\WW$.  
In the first case, $\LL^{(d)}=O_p((d-d_0)^{-1/2})$, so $\|\hat{\boldsymbol{\beta}}^{(d)}_{d_0}-\hat{\boldsymbol{\beta}}_{\lambda_z}\|=O_p((d-d_0)^{-1/2})$; we can also easily obtain the result of \cite{kobak2020optimal} that $\|\hat{\boldsymbol{\beta}}^{(d)}_{d_0}-\hat{\boldsymbol{\beta}}_\lambda\|=o(1)$ a.s., although the rates of convergence give additional information.  For the remaining elements, $\|\hat{\boldsymbol{\beta}}^{(d)}_{-d_0}\| = O_p(1)$ which is not useful.
In the second case,  $\LL_{(n)} = O_p(n^{-1/2})$ and $\|\hat{\boldsymbol{\beta}}_{(n)}-\hat{\boldsymbol{\beta}}_{\lambda}\| = O_p(n^{-1/2})$.
As explained in the Introduction, we argue that our results are more interesting and useful than these.



\subsection{Negative ridging}

In Figure \ref{NegLambdaIdentityCov}, we showed that the optimal value of the ridge penalty in a ridge regression estimator can be negative. Intuitively, this makes sense because including noise in the test set has the same effect as ridge regularisation so, if we have already overshrunk through including too much noise, the optimal action may be to reduce the shrinkage by using a negative ridge parameter. Theorem \ref{thrm_negtv_ridge} below which extends Theorem \ref{thrm_kobak} for sequence II supports this intuition.

Suppose that for sequence II, we set $n_0=n\nu$ with $0 < \nu < 1$.  In Figure \ref{NegLambdaIdentityCov}, this corresponds to making an approximation along the line $n-n_0 =  n_0\nu^{-1}(1-\nu)$ which passes through the origin and has slope $\nu^{-1}(1-\nu)$.  Let $\hat{\be}_{(n),\lambda}$ denote the ridge regression estimator computed along sequence II. Then write $\hat{\boldsymbol{\beta}}_{(n),n\lambda}-\hat{\boldsymbol{\beta}}_{n(1-\nu+\lambda)}= \LL_{(n)}\XX^T\yy$, where the $d_0 \times d_0$ matrix $\LL_{(n)} = (\mathbf{X}^T\mathbf{X}+\WW^T\WW +n\lambda\II_{d_0})^{-1}\{(n-n_0) \mathbf{I}_{d_0}-\WW^T\WW\}\{\XX^T\XX + n(1-\nu+\lambda)\II_{d_0}\}^{-1}$.
We have $\LL_{(n)} = O_p(n^{-3/2})$ and $\|\hat{\boldsymbol{\beta}}_{(n),n\lambda}-\hat{\boldsymbol{\beta}}_{n(1-\nu+\lambda)}\| = O_p(n^{-1/2})$, which proves the following Theorem.
\begin{thrm} \label{thrm_negtv_ridge}
Let $\hat{\boldsymbol{\beta}}_{(n_0),\lambda}=(\mathbf{X}^T\mathbf{X}+\lambda \mathbf{I}_{d_0})^{-1}\mathbf{X}^T\mathbf{y}$ be the ridge regression estimator applied to the data $(\yy,\XX)$. 
Under sequence II with $n_0=n\nu$,  $\|\hat{\boldsymbol{\beta}}_{(n),n\lambda}-\hat{\boldsymbol{\beta}}_{n(1-\nu+\lambda)}\| = O_p(n^{-1/2})$ as $n \rightarrow \infty$.  
\end{thrm}
\noindent The optimal choice of $n(1-\nu+\lambda)$ for the approximating ridge regression estimator $\hat{\boldsymbol{\beta}}_{n(1-\nu+\lambda)}$ is $\|\be_0\|^{-2}d_0\sigma^2$ \citep{nakkiran2020optimal} so the optimal choice of $\lambda$ is $\lambda_{opt} = n^{-1}\|\be_0\|^{-2}d_0\sigma^2 - (1-\nu)$ which can be zero or negative as we showed in Figure \ref{NegLambdaIdentityCov}.

There is an analogous result for the normalized predictors used by \cite{kobak2020optimal} instead of sequence II. Recall that \cite{kobak2020optimal} keep $n_0$ fixed and replace $\WW$ in the training set $\mathcal{D}_{(n)}$ by $\tilde{\WW}=(n-n_0)^{-1/2}\lambda_w^{1/2}\WW$. 
Then we can show that 
$\|\hat{\boldsymbol{\beta}}_{(n),\lambda}-\hat{\boldsymbol{\beta}}_{\lambda+\lambda_w}\| = O_p(n^{-1/2})$.  The optimal ridge parameter in the approximating ridge regression estimator is $\lambda_{opt} = \|\be_0\|^{-2}d_0\sigma^2 - \lambda_w$ which can also be zero or negative, but does not explicitly depend on $n$.  Figure \ref{NegLambdaIdentityCov} shows empirically that the optimal $\lambda$ should depend on $n$, as shown in our theorem.

Theorem \ref{thrm_negtv_ridge} does not hold for sequence I  because the noise predictors are independent (\cite{nakkiran2020optimal}). 



\section{Data Application}
\label{application}
We used publicly available data from the high-density rice array (HDRA, featuring 1,568 accessions and 700,000 SNPs) described by \cite{mccouch2016open} to illustrate our results. The genotype  (allele frequencies) and phenotype (average grain length) data can be downloaded from the Rice Diversity database (\href{}{http://www.ricediversity.org}). We emphasize that our objective is not to reanalyze the dataset but to show that double descent occurs in high-dimensional real-world data.

We reduced the dataset to 332 accessions by only including those from the indica subpopulation with available phenotypic data. We specifically focused on SNPs located on chromosome-3 because a functional SNP in the GS3 gene on rice chromosome-3 was found to be highly associated with grain length in \cite{mccouch2016open}. We removed SNP loci from our analysis if the call rate was less than 95\% and the minor allele frequency (MAF) was less than 5\%.  We highlight the fact that the covariates are non-normal and exhibit high correlation with many instances of perfect correlation due to linkage disequilibrium: strong correlations between SNPs at adjacent loci \citep{ardlie2002patterns, malo2008accommodating}. To reduce the high multicollinearity, we selected every fifth marker to produce a subset of 1,867 markers in our analysis. None of the SNPs has more than 5\% missing values, and we imputed such missing values with the mean of its non-missing values. We standardized each SNP to have zero mean and unit standard deviation and centred the response variable (phenotype data) to have zero mean.  

We randomly split the 332 accessions, allocating 300 to the training and 32 to the test sets. Starting from an intercept-only regression model, we gradually added genetic markers in their natural sequence to assess the trends in both training and test errors (as outlined in steps 9 to 12 of Algorithm S1). This process was iterated 100 times, each with new random data split into training and testing sets. 

\begin{figure}[t!]
\begin{center}$
\begin{array}{c}
\includegraphics[width=6cm, height=10cm, keepaspectratio]{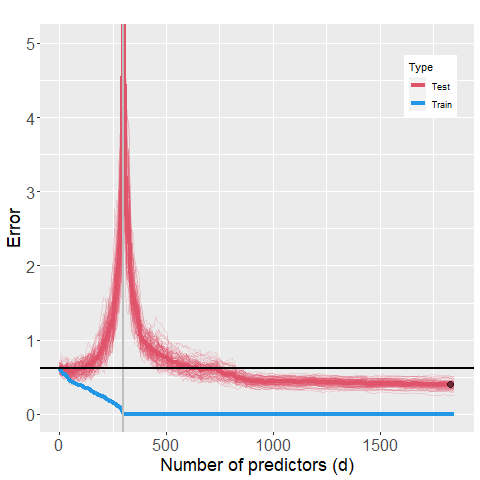}
\end{array}$
  \caption{Root mean squared error (RMSE) from regression analysis predicting average grain length in Indica rice using selected genetic markers on chromosome-3. }\label{GWAS_error_plot}
\end{center}
\end{figure}

The test error in Figure \ref{GWAS_error_plot} shows double descent, reaching its minimum in the over-parameterized regime. The noticeable decrease in test error just after $d=800$ highlights a SNP previously identified as significantly associated with grain length.

\section{Discussion}
\label{sec:conc}

Most datasets are assumed to include noise predictors; this is encapsulated in the widely used assumption of sparsity and underlies the use of model selection methods which try to identify the important predictors. Datasets may also include noise observations.  Our study of the effect of including noise predictors and noise observations when fitting linear regression models show that both forms of noise perform shrinkage in a similar way to ridge regression estimators.  We develop the consequences of this shrinkage due to noise to provide simple, intuitive explanations for double descent (in terms of the asymptotes of the test error and the ``condition number anomaly'') and negative ridging (in terms of ``double shrinkage''). These explanations also apply to shrinkage due to measurement error \citep{webb1994functional, bishop1995training, hastie2022surprises}.

Moreover, the fact that double descent occurs both in the overparameterized regime (when fitting noise predictors) and in the underparameterized regime (when fitting noise observations) shows that it is not directly driven by overparameterization (which occurs only in the first case), but rather is driven by the shrinkage (which occurs in both cases) due to fitting noise.  

Our results raise the question of whether, under sparsity, it is better to select a model (often in the underparameterized regime) or to fit a (regularized) model in the overparameterized regime. In so far as lasso regression \citep{tibshirani1996regression} represents the first approach, the question can be formulated as: Should we prefer lasso or ridge regression?  It is possible that the lasso's ability to set some components of $\hat{\boldsymbol{\beta}}$ to zero may be less of an advantage than has hitherto been thought. Finally, in addition to having potentially better performance, the fact that shrinkage provides a way to interpret overparameterized models may help to reduce or even overcome the claimed interpretability advantages of simple models.  Addressing these questions and extending our results to other models, such as neural networks, are potential avenues for future research.

\bigskip
\begin{center}
{\large\bf SUPPLEMENTARY MATERIAL}
\end{center}

\begin{description}

\item[Supplementary material:] The supplementary materials include the detailed algorithms to produce the figures, additional figures, technical lemmas along with their proofs, as well as the proofs of the theorems presented in the main text. (pdf file)

\item[R code:] The zipped folder contains R scripts for the exact replication of all analyses and figures from the main text and supplementary materials. (zipped file)

\item[High-density rice array dataset:] 
The dataset used to illustrate the method is publicly available at \href{}{http://www.ricediversity.org}. The accompanying zipped file includes the R code to download the data and reproduce the analysis presented in section \ref{application}.
\end{description}

\bibliographystyle{agsm}


\end{document}